\documentclass[12pt]{article}
\usepackage{amsmath}
\usepackage{graphicx,psfrag,epsf}
\usepackage{enumerate}
\usepackage{natbib}
\usepackage{url} 
\usepackage[utf8]{inputenc}
\usepackage{amsmath,amsthm}
\usepackage{amsfonts,graphicx}
\usepackage{amssymb,hyperref}
\usepackage{color,url}
\usepackage{natbib}
\usepackage{todonotes}
\usepackage{soul}
\newtheorem{prop}{Proposition}
\newtheorem{thm}{Theorem}
\newtheorem{lema}{Lemma}
\newtheorem{remark}{Remark}

\newcommand{\eps}{\varepsilon}

\newcommand{\ind}{\mathbb{I}}

\newcommand{\Haus}{\mathbb{H}}

\newcommand{\blind}{0}

\begin{document}

\def\spacingset#1{\renewcommand{\baselinestretch}%
{#1}\small\normalsize} \spacingset{1}


\if0\blind
{
  \title{\bf On visual distances for spectrum-type functional data}
  \author{Alejandro Cholaquidis$^a$, Antonio Cuevas$^b$ and Ricardo Fraiman$^a$.  
   \\ \\
   $^a$ Universidad de la Rep\'ublica.\\
 $^b$  Universidad Aut\'onoma de Madrid.} 
  \maketitle
} \fi

\if1\blind
{
  \bigskip
  \bigskip
  \bigskip
  \begin{center}
    {\LARGE\bf On visual distances for spectrum-type functional data}
\end{center}
  \medskip
} \fi

\bigskip

\begin{abstract}
	A functional distance ${\mathbb H}$, based on the Hausdorff metric between the function hypographs, is proposed for the space ${\mathcal E}$ of non-negative real upper semicontinuous functions on a compact interval. The main goal of the paper is to show that the space $({\mathcal E},{\mathbb H})$ is particularly suitable in some statistical problems with functional data which involve functions with very wiggly graphs and narrow, sharp peaks. A typical example is given  by spectrograms, either obtained by magnetic resonance or by mass spectrometry. On the theoretical side, we show that $({\mathcal E},{\mathbb H})$ is a complete, separable locally compact space and that the ${\mathbb H}$-convergence of a sequence of functions implies the convergence of the respective maximum values of these functions. The probabilistic and statistical implications of these results are discussed in particular, regarding the consistency of $k$-NN classifiers for supervised classification problems with functional data in ${\mathbb H}$. On the practical side, we provide the results of a small simulation study and check also the performance of our method in two real data problems of supervised classification involving  mass spectra.
\end{abstract}
\noindent\it Key words\rm: Supervised classification, functional data analysis, Hausdorff metric.
\section{Introduction: the choice of a suitable functional distance}

The statistical analysis of problems where the sample data are functions is often called Functional Data Analysis (FDA). This is a relatively new statistical field which involves several specific challenges, most of them are associated with the infinite-dimensional nature of the data.

We are concerned here with one of these specific challenges, namely, the choice of a suitable distance criterion between the data. In what follows, unless otherwise stated, we will consider problems where the sample data are real functions $x(t),\ t\in[0,1]$.

Not surprisingly, a considerable part of the current FDA theory has been developed assuming that the data functions belong to the space $L^2[0,1]$, that is, the distance
between two data $x_1$ and $x_2$ is given by $d_2(x_1,x_2)=(\int_0^1(x_1(t)-x_2(t))^2dt)^{1/2}$. This distance presents obvious advantages, derived from the fact that $L^2[0,1]$ is a Hilbert space. Thus, some extremely important tools, as the existence of orthogonal bases (and the corresponding expansions for the data in orthogonal series) are available in $L^2[0,1]$. As a useful by-product, some crucial methodologies, such as Principal Components Analysis or  Linear Regression (and even Partial Least Squares), can be partially adapted to the functional setting.

Another widely used metric is associated with the supremum norm $\Vert x\Vert_\infty=\sup_t|x(t)|$, which is well-\-defined in the space ${\mathcal C}[0,1]$ of real continuous functions $x:[0,1]\rightarrow {\mathbb R}$; 
thus the metric is $d_\infty(x_1,x_2)=\sup_t|x_1(t)-x_2(t)|$ for $x_1,x_2\in{\mathcal C}[0,1]$. Although the Hilbert structure is lost here, the advantages of the supremum metric  are also well-known: first, $d_\infty$ is easy to interpret in terms of vertical distance between the functions. Second, the structure of the space  of probability measures on $({\mathcal C}[0,1],\Vert\cdot\Vert_\infty)$ is also well understood, and carefully analyzed, for example, in the classical book by \citet{bil68}.

For general accounts on the FDA theory we refer to the books by \citet{bos00}, \citet{bos07}, \citet{ram02,ram05}, \cite{fer06}, \citet{hor12} and the recent survey paper by \citet{cue14}.

\subsection{Our proposal: its practical motivation}

In what follows we analyze, from both the theoretical and practical point of view, a metric between functions especially aimed at capturing the ``visual distance'' between the graphs. This metric will be particularly suitable in FDA problems where the data are functions with wiggly graphs showing very sharp peaks. In those situations the classical metrics ($d_2$ or $d_\infty$)  could be unsuccessful  in capturing a ``practically meaningful'' notion of distance between the graphs. For example, a small lateral shift in a very sharp peak (perhaps due to a registration error) could lead to an enormous $d_\infty$-distance. Likewise, if two graphs differ in just one such narrow peak, the $d_2$-distance between them might be very small, which could be unsuitable in many cases.

\begin{figure}[h]
	\centerline{\includegraphics[scale=0.33]{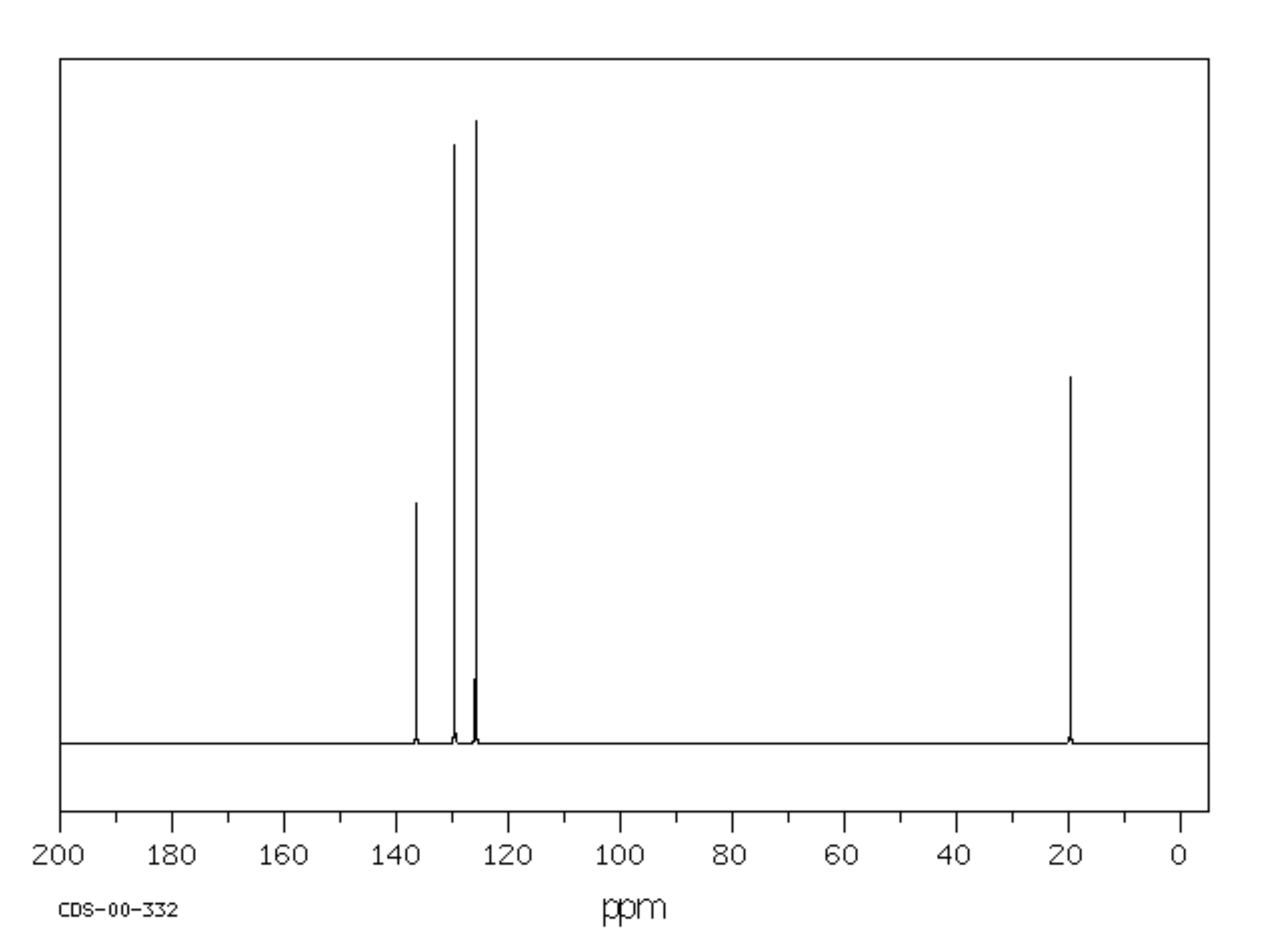}\includegraphics[scale=0.33]{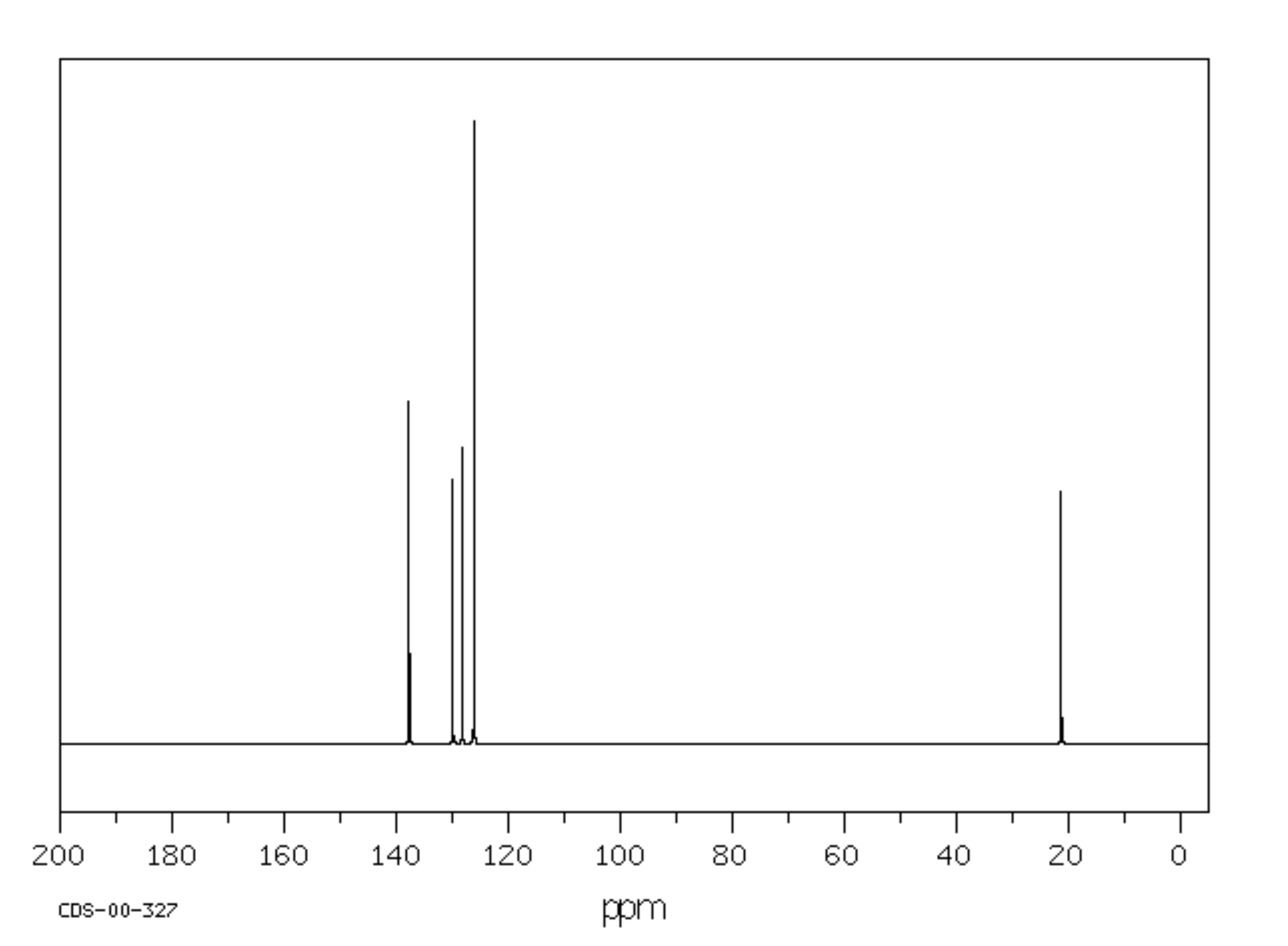}}
	\caption{\footnotesize $^{13}$C-NMR spectra of the o-xylene (left) and m-xylene.} 
	\label{xylene}
\end{figure}

The spectrograms, either obtained from magnetic resonance ($^1$H-NMR or $^{13} $C-NMR) or by mass spectrometry, provide a good example of such situations. Just as an example to motivate our point, let us consider the $^{13} $C-NMR spectrum of a compound, namely the o-xylene ($C_8H_{10}$); see Figure \ref{xylene}, left. It shows the typical spiky pattern, with sharp and narrow peaks, strongly localized (we will consider below other examples of  much more complex organic compounds were the peaks are present but not all the information is concentrated around them). The peaks   in this spectrum are located at the points 136.42, 129.63, 125.85, 19.66 ppm.  This information has been obtained from the data base \url{http://sdbs.db.aist.go.jp}, (National Institute of Advanced Industrial Science and Technology, date of access August 23, 2015). 
Now, we might want to consider the $^{13}$C-NMR spectrum of another closely related compound, the m-xylene, an isomer of the previous one: see the right panel of Figure \ref{xylene}. Although the general aspect of both spectra is very similar, there are clearly some differences. In the case of m-xylene, the peaks are located at 137.74, 129.96, 128.21, 126.13, 21.31 ppm.  A  ``reasonable'' metric defined to measure the distance between these graphics should provide a small value (thus reflecting their close affinity), by taking into account their ``visual'' proximity, that is, the distance between the graphics in all directions (not only in the vertical one). Moreover, for this type of graphics, we would also like to detect the presence of additional
very narrow peaks  (far away from the others), contributing a small area but carrying a relevant information on the compound. The $d_2$ distance does not seem useful for such purpose. See also the Figure \ref{disthaus} below and the discussion following definition (\ref{HH}).

As explained in depth by \citet{ba07}, in order to reach meaningful conclusions, handling of spectrum data needs a crucial pre-\-processing stage. This typically includes, among others, the following steps: \textit{remove random noise}, \textit{normalization}, \textit{peaks detection} (to identify locations on the scale that correspond to specific molecules) and \textit{peak matching} (to match peaks in different samples, that correspond to the same peak). For this purpose, there is an increasing amount of software available. In particular, several packages can be downloaded from the web page of the software \tt R \rm (\url{http://www.r-project.org/}) in order to deal with spectrum-type data; for example, \tt MALDIquant\rm,   \tt readMzXmlData\/ \rm and \tt aLFQ\rm. This paper could be seen as a further suggestion in this line of research.

The rest of this paper is organized as follows: the study of the proposed visual metric (including the definition, computation and topological properties of the distance) is considered in Section \ref{thedistance}. In Section \ref{classification} we focus on some theoretical aspects of the use of this metric in the supervised classification problem. A small simulation study is provided in Section \ref{simul}.  Two real data examples of mass spectra classification are considered in Section \ref{real}.  Finally, some concluding remarks are given in Section \ref{conclu}. The proofs are given in an appendix.

\section{A visual, Hausdorff-based distance for non\--negative functions}
\label{thedistance}

The starting point is the standard definition (see, e.g., \citet{roc09}, p. 117)  of the Hausdorff (or Pompeiu-Hausdorff) distance between two compact non-empty sets, $A,C\subset{\mathbb R}^d$:
\begin{eqnarray*}
	d_H(A,C)&=&\inf\{\epsilon>0:\, A\subset B(C,\epsilon),\, C\subset B(A,\epsilon)\}\\ \noindent
	& =&\max \bigg\{ \max_{a\in A}d(a,C),\ \max_{c\in C} d(c,A)\bigg\},\label{H}
\end{eqnarray*}
where $B(A,\epsilon)$ denotes the $\epsilon$-parallel set $B(A,\epsilon)=\cup_{x\in A}B(x,\epsilon)$ and $B(x,\epsilon)$ denotes (with a slight notational abuse) the closed ball centered at $x$ with radius $\epsilon$; the open ball will be denoted $\mathring{B}(x,\epsilon)$. Also, $d(a,C)=\inf\{\Vert a-c\Vert:\, c\in C\}$.

Unlike other notions of proximity between sets, $d_H$ is a true metric (i.e. it has the properties of identifiability, symmetry and triangular inequality) in the class of compact non-empty sets. The Hausdorff distance has been extensively used in different problems of image analysis (especially in pattern recognition), which appear in the literature under different names (shape comparison, object matching, etc.). The strong intuitive motivation behind the definition of $d_H$ has motivated the study of other variants of the same idea as well as other closely associated notions. Some  references are \cite{hut93}, \cite{dub94}, \cite{sim99}.

However, our aim here is rather to use the Hausdorff metric as a tool for defining distance between functions, very much in the spirit of some ideas in 
approximation theory; see, for example, \cite{sen90}. 

The basic idea behind the metric we are going to consider is quite simple: given two non-negative functions $f$ and $g$,
defined on $[0,1]$, the distance between $f$ and $g$ is measured in terms of the Hausdorff metric between the corresponding hypographs. However, we must take care of some technical aspects in order to properly establish this definition.

Let us recall that a function $f:[0,1]\rightarrow {\mathbb R}$ is said to be \it upper semicontinuous at $x_0$\/ \rm if $\limsup_{x\rightarrow x_0}f(x)\leq f(x_0)$. A
function  $f$ is said to be upper semicontinuous (USC) if it fulfils the above condition at every point $x_0$.

Given a non-negative function $f$ defined on $[0,1]$, the \it hypograph of $f$ \/ \rm is the set
$$
H_f=\big\{ (x,y)\subset {\mathbb R}^2: x\in [0,1], 0 \leq  y \leq f(x)\big\}.
$$
Denote by $\mathcal{E}$ the space of non-negative USC functions defined on $[0,1]$. The following proposition, whose proof can be found in \cite{nat60}, establishes
some useful properties of USC functions.

\begin{prop} \label{propusc} Let $f:[0,1]\rightarrow \mathbb{R}$ be a non-negative USC function.
	\begin{itemize}
		\item[1)] Let $K \subset [0,1]$ be  a compact set. Then, there exists $z\in K$ such that $\sup_{x\in K} f(x)=f(z)$.
		\item[2)] The hypograph $H_f$ is compact.
	\end{itemize}
\end{prop}

We are now ready to define our visual metric: for $f$, $g\in{\mathcal E}$ we define
\begin{equation}
	{\mathbb H}(f,g)=d_{H}(H_f,H_g)\label{HH}.
\end{equation}
It is easily seen that this is a true metric in ${\mathcal E}$. In particular, if $d_{H}(H_f,H_g)=0$, the USC assumption
guarantees that we must have $f=g$. 

Let us denote by $({\mathcal E},{\mathbb H})$ the space of USC non-negative functions endowed with the metric (\ref{HH}).

\begin{figure}[h]
	\centerline{ \includegraphics[scale=0.24]{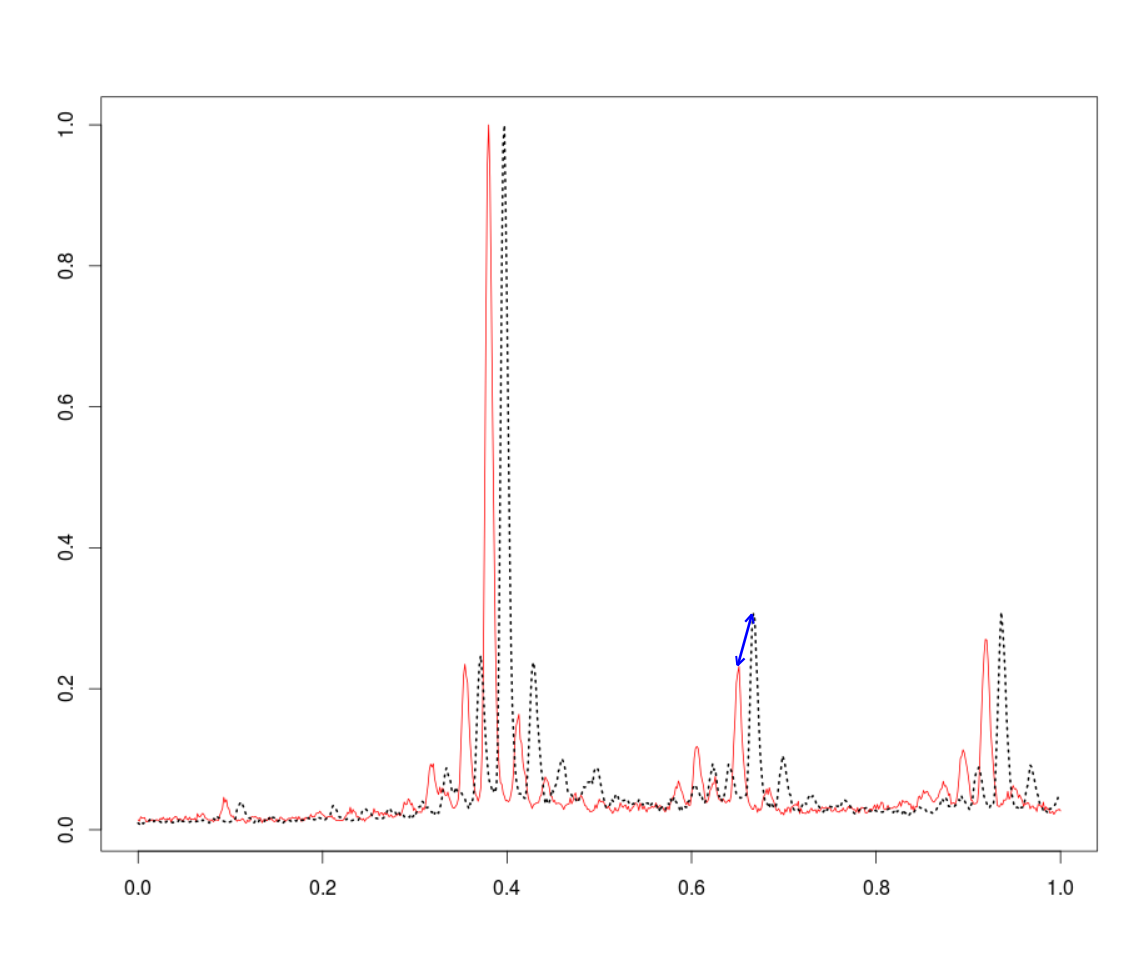}}
	\caption{$\mathbb{H}(f_1,f_2)$ corresponds to the length of the double arrow joining both curves.}
	\label{disthaus}
\end{figure}

Figure \ref{disthaus} aims at illustrating  the heuristic meaning of  the ${\mathbb H}$-metric between two spectrum-type curves $f_1$ (solid line) and $f_2$ (dotted line). In this case, the distance $\mathbb{H}(f_1,f_2)=0.077$ is the length of the double arrow joining two close peaks in the curves. The corresponding  values for the ``classical" distances are  $d_\infty(f_1,f_2)=0.952$, $d_2(f_1,f_2)=0.119$. So, ${\mathbb H}$ succeeds in reflecting the visual proximity between both curves.

\begin{remark}
	In order to gain some additional insight on the meaning of the distance ${\mathbb H}$ and their relations with other usual metrics, let us note that
	\begin{itemize}
		\item[(a)] Convergence in ${\mathbb H}$ does not imply pointwise convergence. Consider  $f(x)={\mathbb I}_{[0,1]}(x)$ and
		the sequence  $f_n(x)=nx$ if $x<1/n$ and $f_n(x)=1$
		if  $x\in [1/n,1]$. It is clear that  $\Haus(f_n,f)\rightarrow 0$ but $f_n(0)\nrightarrow f(0)$. The reciprocal implication is not true either. Take $f(x)=0$ $\forall x\in [0,1]$
		and $f_n(x)=x^n(1-x^{n})$. We have $f_n(x)\to 0$ for all $x$
		but  $\Haus(f_n,f_0)\rightarrow 0$ where $f_0(1)=1/4$  and $f_0(x)=0$ for $x\in [0,1)$.
		\item[(b)] Convergence in $L^p$ does not imply convergence in $\Haus$: consider $f_n(x)=\ind_{[0,1/n]}$ and $f(x)=0$ $\forall x\in [0,1]$. The reciprocal is also false: take
		$I^n_i(x)=\ind_{\left[\frac{i-1}{2^n},\frac{i}{2^n}\right]}(x)\quad i=1,\dots,2^n$,
		and the functions $$f_n(x)=\sum_{k=1}^{2^n-2}I^n_{2k}(x)\quad \text{and}\quad f(x)=\ind_{[0,1]}(x),$$
		it is clear that  $\Haus(f,f_n)=1/2^{n+1}$, however, for all $p$, $\int^1_0|f_n(x)-f(x)|^pdx=1/2$.
		
		\item[(c)] A natural question is: why to use USC functions? Since in many applications, the functional data appear as continuous functions, one might think that we might restrict our discussion to the continuous case. However there are, at least, good reasons for considering USC functions: first, in some practical examples (though not, typically, in the case of spectra) one has to deal with non-continuous samples and, especially, with jump discontinuities. Second, we need upper semicontinuity to get a complete space (which in turn is essential for a more convenient mathematical handling). Take for example the sequence of functions $f_n(x)=x^n$ for $x\in[0,1]$. This is a Cauchy sequence in our space $({\mathcal E},{\mathbb H})$, but clearly does not converge to any continuous function on [0,1]. Hence, we need to enlarge the space to include USC if we want to get completeness. 
		
	\end{itemize}

\end{remark}

\subsection{Computational aspects}

As mentioned above, the Hausdorff distance has some applications in image processing. Hence its numerical calculation has
motivated some interest in the literature.  See \cite{nu11} and \cite{alt03}, just to mention a couple of recent references.
The \tt Matlab \rm function  \tt HausdorffDist \rm computes the Hausdorff distance between two finite sets of points in $\mathbb{R}^2$.
We are concerned here with the particular case in which the sets are the hypographs of functions, especially when these functions are specified only by their values in a given grid of $[0,1]$. Our aim here is to approximate the value of ${\mathbb H}(f,g)$ in such cases, which arise very often in practical applications.  First let us observe that, given two functions $f$ and $g$
in $\mathcal{E}$, from the definitions of ${\mathbb H}$ and $d_H$ we have: $\Haus(f,g)\leq d_H(\partial H_f, \partial
H_g)$. However,
the boundaries $\partial H_f$ and $\partial H_g$ do play a relevant role in the calculation
of ${\mathbb H}(f,g)$. In fact, the following proposition shows that we can restrict the calculation to appropriate subsets of these boundaries.

\begin{prop} \label{prophaus}
	Let $f,g\in \mathcal{E}$, then
	\footnotesize{$$\Haus(f,g)=\max \left\{\sup_{\left\{\substack{x=(x_1,x_2)\in \partial H_g \\ g(x_1)\geq f(x_1)}\right\}}d\big(x,\partial H_f\big),\ \sup_{\left\{\substack{y=(y_1,y_2)\in \partial H_f \\ f(y_1)> g(y_1)}\right\}} d\big(y,\partial H_g\big)\right\}.$$}
\end{prop}

The proofs of all results are given in the Appendix. In particular, the proof of Proposition \ref{prophaus}, will require Lemma \ref{lemhaus} whose proof can also be found in the Appendix. 

\ 

\noindent \textit{An algorithm to compute $\Haus$}\\

If we have $\alpha_n^ f=(t_1,f(t_1)),\dots,(t_n,f(t_n))$ and $\alpha_n^g=(t_1,g(t_1)),\dots,(t_n,g(t_n))$, and if we assume that $f$ and $g$ are  continuous, then Proposition \ref{prophaus}, together with Proposition \ref{propdh} in the Appendix, gives us a simple algorithm of order $n^2$ to approximate $\Haus(f,g)$. Just compute
$$\tilde{\Haus}(f,g)=\max \bigg\{\max_{\{i:g(t_i)> f(t_i)\}}d\big((t_i,g(t_i)),\alpha_n^f\big),\ \max_{\{i:f(t_i)> g(t_i)\}} d\big((t_i,f(t_i)),\alpha_n^g\big)\bigg\},$$
where, if $\{i:g(t_i)>f(t_i)\}=\emptyset$ then $\max_{\{i:g(t_i)> f(t_i)\}}d\big((t_i,g(t_i)),\alpha_n^f\big)=0$ and if $\{i:f(t_i)>g(t_i)\}=\emptyset$,  then $\max_{\{i:f(t_i)> g(t_i)\}} d\big((t_i,f(t_i)),\alpha_n^g\big)=0$. From  Lemma \ref{lemhaus}, if $\max_i{|t_{i+1}-t_i|}\rightarrow 0$ then $\tilde{\Haus}(f,g)\rightarrow \Haus(f,g)$.

\subsection{Some related literature}

The distance ${\mathbb H}$ has been considered in \cite{cue98}
in the context of density estimation: in particular, convergence rates are obtained, under some smoothness conditions,  for ${\mathbb H}(\hat f_n,f)$, where $\hat f_n$ denotes a sequence of kernel density estimators of the density $f$.

Different versions of the same idea are considered in \cite{roc09}, p. 282. They are defined in terms of epigraphs (rather than hypographs) and are therefore applied to lower semicontinuous (rather than upper semicontinuous) functions. Some relevant applications are given in the framework of optimization theory to give bounds for approximately optimal solutions of convex lower semicontinuous functions.

Another related approach to the idea of defining the distance between two functions in terms of the distance between their graphs is considered in \cite{sen90} for the so-called \it segment functions\rm. \cite{hol92} extends these ideas to the setting of multifunctions.

\subsection{Topological properties of $({\mathcal E},{\mathbb H})$}


The metric space is particularly ``well-behaving'' in some important aspects that are summarized next.

\begin{thm}\label{topo}
	(a) The space $({\mathcal E},{\mathbb H})$ is complete and separable.  Also, any bounded and closed  set in $({\mathcal E},{\mathbb H})$ is compact. In particular, $({\mathcal E},{\mathbb H})$ is locally compact.
	
	(b) Let $f_n, f\in {\mathcal E}$ such that $\Haus(f_n,f)\rightarrow 0$ then;
	$$\max_{x\in [0,1]} f(x)= \lim_{n\rightarrow +\infty} \max_{x \in [0,1]} f_n(x).$$
\end{thm}

\

The proof of this result is given in the Appendix. Let us now briefly comment on the meaning and usefulness of these properties.

\begin{itemize}
	\item[(i)] Among the three properties established in Theorem \ref{topo} (a), completeness is perhaps the most basic one. It is essential to study convergence of sequences or series in $({\mathcal E},{\mathbb H})$ by just looking at the corresponding Cauchy property. This property is also required in the proof of some key results as Banach's fixed point theorem for contraction mappings; see \cite{granas:03}
	
	\item[(ii)] Separability is a most crucial property in a metric space in order to define on it well-behaving probability measures. A nice discussion on this topic can be found in \cite{led91}, pp. 38-39. Although this discussion applies, in principle, to Banach spaces, the main arguments can be also translated to a metric space. For example, separability is required to ensure that a probability measure $P$ defined on $({\mathcal E},{\mathbb H})$ is \it tight\rm, in the sense that for all $\epsilon>0$ there exists a compact set $K\subset {\mathcal E}$ such that $P(K)>1-\epsilon$. This is a far-reaching property, that can be found in the basis of many standard probability calculations. Thus, the separability property allows us to express any ${\mathcal E}$-valued random element as a limit of a sequence of simple (finite-valued) random elements. Also, separability is needed to guarantee a proper behaviour of product measurable structures: in particular, the Borel $\sigma$-algebra of the product space is the product of the individual Borel $\sigma$-algebras of the factors; see Proposition 1.5 in \cite{fol99}.
	Also, let us recall that separability of a metric space is equivalent to the property that this space is second-countable (e.g., \cite{fol99}, pp. 116-118), which is important in many probability arguments: for example, to show that any probability measure in a locally compact metric space is a Radon measure, see \cite{fol99}, Th. 7.8.
	
	Finally, separability is also required for the consistency theorem for $k$-NN classification rules mentioned in Subsection \ref{acfd}.
	\item[(iii)] As for local compactness, let us recall that, in the case of Banach spaces, this property is equivalent to the finite-dimensionality of the space. In our case, we don't have a vector structure, so that we only have a metric space (not a normed one). However, the local compactness allows us to use some ``natural'' properties that we often use in the finite-dimensional spaces.
	For example, to show that any real integrable function defined on ${\mathcal E}$
	can be approximated by a sequence of continuous compact-supported functions (see \cite{fol99}, Proposition 7.9). An application of this can be found in Section \ref{ppp}.
\end{itemize}

\begin{remark}
	Let us observe that the local compactness does not hold for $(\mathcal{E},\Vert\cdot\Vert_\infty)$. In order to see this, observe that for every $\epsilon>0$,
	the sequence $f_n(x)=\epsilon x^n$
	is included in the ball (with the distance $\|\cdot\|_\infty$) centered at the null function, of radius $\epsilon$. However, this sequence does not have any convergent
	subsequence; indeed, the only possible limit would be the function $f_0(x)=0$ for $x\in[0,1)$, $f_0(1)=\epsilon$, but  $\|f_n-f_0\|_\infty=\epsilon$ for all $n$.
\end{remark}

\section{Applications to classification of functional data}\label{classification}\label{acfd}
We will briefly consider here some theoretical aspects of the supervised classification problem, focusing especially  on the case of $k$-NN (nearest neighbors) classifiers.

\

\noindent \it The (functional) supervised classification problem\rm

We focus on the problem of  supervised classification with functional data; see e.g., \cite{bai11a} for an overview. More precisely, we are concerned with statistical problems for which the available data consist of an iid ``training sample'' $D_n=((X_1,Y_1),\ldots,(X_n,Y_n))$. The $X_i=\{X_i(t):\,t\in[0,1]\}$ are independent trajectories, belonging to a function space ${\mathcal X}$, drawn from a stochastic process $X:=\{X(t):\, t\in[0,1]\}$ which can be observed from two probability distributions, $P_0$ and $P_1$ (often referred to as ``populations'' in statistical language). The $Y_i$ are binary random variables indicating the membership of the trajectory $X_i$ to $P_0$ or $P_1$, that is, the population from which the observation $X_i$ has been drawn. It is assumed that the conditional distributions of $X|Y=i$, $i=0,1$ (that is, $P_0$ and $P_1$) are different.

\

\noindent \it $k$-NN classifiers: why to use them in the functional setting\rm

In a model of this type, the aim is typically to classify (either in $P_0$ or in $P_1$) a new observation $X$, for which the corresponding value of $Y$ is unknown. A classification rule (or classifier) is a measurable function $g:{\mathcal X}\rightarrow \{0,1\}$ defined on the space ${\mathcal X}$ of trajectories. Usually, the classification rules are constructed using the information provided by the training sample data $(X_i,Y_i)$.

In this work we will limit ourselves to use the $k$-NN classifiers: an observation $x$ is classified into $P_0$ if the majority among the $k$  observations $X_i$ (in the training sample) closest to $x$, fulfils $Y_i=0$; ties are randomly broken. Of course, ``closest'' refer to some metric defined in the space ${\mathcal X}$ on which the $X_i$ take values: each metric leads to a different $k$-NN classification rule. In the functional infinite dimensional case, the choice of this metric is particularly relevant.
The values $k=k_n$ are the smoothing parameters, similar to others which appear in non-parametric procedures: see \cite{dev96} for background. As we will see below, they must fulfil some minimal conditions regarding the speed of convergence to infinity.
Of course, the choice of $k$ for any specific sample size $n$ can have some influence on the performance of the $k$-NN
classifier. However, as we will see in
Section \ref{real}, the choice of the metric in the ``feature space'' (where $X$ takes values) can be even more important.

The reasons for choosing $k$-NN classifiers can be summarized in the following terms: simplicity, ease of interpretability, good general performance and generality. Indeed, $k$-NN is a sort of all-purposes ``benchmark procedure'', not so easy to beat in practice. The available experience (see \citet[2011a]{bai11}, \cite{gal14} and references therein) suggests that,  $k$-NN classifiers tend to show a stable performance, not far from the best method found in every specific problem. Moreover, they have a sound intuitive basis, so they are easily interpretable in all cases (unlike other classification methods) and they can be used in very general settings, when $X$ takes values in any metric space.
We now consider some theoretical issues regarding consistency
of $k$-NN classifiers  in the framework of our space $({\mathcal E},{\mathbb H})$.

\

\noindent \it The notion of consistency\rm

Let us denote by $g_n$ a sequence of $k$-NN classifiers defined in the usual way, as indicated at the end of the previous subsection. We will say that this sequence is \it weakly consistent\/ \rm (see, e.g., \cite{dev96} for more details) if the misclassification probability $L_n={\mathbb P}(g_n(X)\neq Y|D_n)$  converges (in probability, as $n\to\infty$) to the optimal value $L^*={\mathbb P}(g^*(X)\neq Y)$, which corresponds to the optimal rule
$g^*(x)={\mathbb I}_{\{\eta(x)>1/2\}}$, where
$\eta(x)={\mathbb E}(Y|X=x)={\mathbb P}(Y=1|X=x)$.
It is readily seen that the weak consistency condition is equivalent to
\begin{equation*}
	{\mathbb E}(L_n)\longrightarrow {\mathbb E}(L^*).\label{wc}
\end{equation*}

In the finite dimensional case, that is when random variable $X$ takes values in ${\mathbb R}^d$,  it is well-known from a classical result due to \cite{sto77}, that any sequence of $k$-NN classifiers is (weakly) consistent provided that $k\to\infty$ and $k/n\to 0$. This result is \it universal\rm, in the sense that it does not impose any condition of the distribution of the random pair $(X,Y)$.

\

\noindent \it The infinite-dimensional case. The Besicovitch condition\rm

Let $(X,Y)$ be the random element generating the data in a supervised functional classification problem, where $X$ is ${\mathcal E}$-valued  and $Y$ takes values in $\{0,1\}$. Denote by $\mu$ the distribution of $X$, $\mu(E)={\mathbb P}(X\in E)$.

It is natural to ask whether the above mentioned universal consistency of the finite-dimensional $k$-NN classifiers still holds for the functional (infinite-dimensional) case. The answer is negative. There is, however, an additional technical condition which (together with $k\to\infty$, $k/n\to 0$), ensures weak consistency for the $k$-NN functional classifiers. While this condition is not in general trivial to check, it always holds whenever the regression function $\eta(x)={\mathbb E}(Y|X=x)$ is continuous. The corresponding theory has been first developed by \cite{cer06}. In particular, the mentioned sufficient condition for consistency established by these authors is the following differentiability-type assumption (on the distribution of $(X,Y)$), called \it Besicovitch  condition\rm:

\begin{equation} \label{besicovitch2}
	\lim_{\delta\rightarrow 0} \frac{1}{\mu\big(B(X,\delta)\big)}
	\int_{B(X,\delta)} |\eta(X)-\eta(x)|d\mu(x)=0,\ \mbox{in probability}. \
\end{equation}
Here, $B(X,\delta)$ denotes the closed $\delta$-ball centered at $X=X(t)$ in the space of trajectories of the process $X=X(t)$. A weaker, slightly simpler version of this property, almost identical to the conclusion of Lebesgue differentiation theorem, would be as follows,
\begin{equation} \label{besicovitch1}
	\lim_{\delta\rightarrow 0} \frac{1}{\mu\big(B(X,\delta)\big)}
	\int_{B(X,\delta)} \eta(x) d\mu(x)=\eta(X),\ \mbox{in probability}. \
\end{equation}

Conditions (\ref{besicovitch2}) and (\ref{besicovitch1}) are clearly reminiscent of the conclusion of the classical Lebesgue Differentiation Theorem (see \citep[p. 98]{fol99}). Clearly (\ref{besicovitch2}) implies (\ref{besicovitch1}). It can be also seen that the $\mu$-a.s. continuity of $\eta$ is a sufficient condition for (\ref{besicovitch2}).

As mentioned above, \citet[Th. 2]{cer06} have proved that condition (\ref{besicovitch2}) together with $k\to\infty$ and $k/n\to 0$, ensures the weak consistency of a sequence of
$k$-NN classifiers when $X$ takes values in a separable metric space. On the other hand, \cite{abr06} have used  (\ref{besicovitch1}) as a sufficient condition for the consistency of kernel classification rules. They also need some supplementary conditions on the sequence $h=h_n$ of smoothing parameters and the space ${\mathcal E}$: they require that the existence of a sequence of non-decreasing totally bounded subsets, $\mathcal{F}_k\subset\mathcal{E}$, such that $\mu(\cup_k \mathcal{F}_k)=1$ and a condition that relates the bandwidth $h$ with the metric entropy of the subsets $\mathcal{F}_k$.

The following result shows that, in our case, the consistency holds for a class of ``regular'' distributions which is dense in the space of all distributions.  In other words, the result shows that the assumption of continuity for the regression function $\eta(x)$ (which guarantees consistency for $k$-NN classifiers) is in fact not very restrictive, as any possible distribution for $(X,Y)$ may be arbitrarily approximated by another one which fulfils this continuity condition.

\begin{prop}\label{ppp}
	Let us consider a binary supervised classification problem based on observations from $(X,Y)$, where $X$ is ${\mathcal E}$-valued and $Y$ is the binary variable indicating the class (0 or 1). Let $g_n$ be a sequence of $k$-NN classifiers such that $k\to \infty$ and $k/n\to 0$.
	
	Whatever the distribution $Q$ of $(X,Y)$ there is another distribution $P$, arbitrarily close to $Q$ in the weak topology, under which the regression function $\eta(x)={\mathbb P}(Y=1|X=x)$ is continuous with compact support and the sequence $g_n$ is weakly consistent.

\end{prop}

\section{Some simulations}\label{simul}

A simulation experiment has been carried out to illustrate a simple situation in which our ``visual distance'' could be especially suitable. The underlying model is very simple: the functional data are just ``corrupted'' trajectories of the absolute value of a Brownian Bridge (absBB) on $[0,1]$. In the population $P_0$ the absBB trajectories are perturbed by just adding to them a spiky function identically null on $[0,1]$ except for a triangular peak with basis 0.04 and height 1, whose center is randomly chosen on the interval $[0,a_1]$ for some $a_1\leq 1/2$. The trajectories from $P_1$ are similarly constructed except that the center of the noise peak is randomly selected on $[a_2,1]$ for some $1/2\leq a_2\leq 1$.

We have performed this experiment for two choices of $(a_1,a_2)$. The first case (Model 1, Table \ref{tablesim1} left) corresponds to the choice $a_1=a_2=1/2$. In the second one (Model 2, Table \ref{tablesim1} right), we have taken $a_1=1/3$, $a_2=2/3$.

\begin{figure}[h]
	\centerline{\includegraphics[scale=0.4]{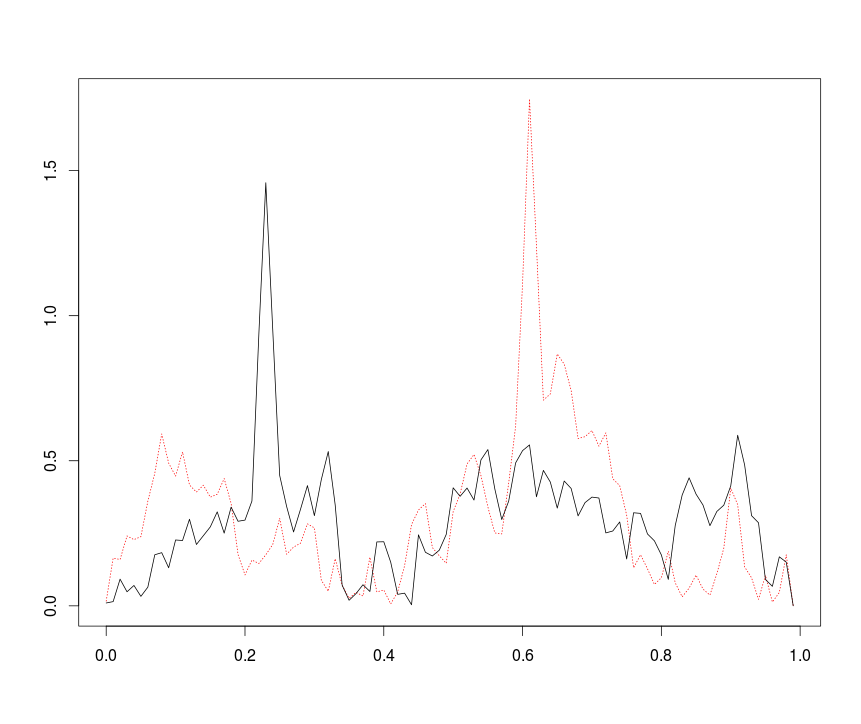}}
	\caption{\footnotesize Two random trajectories drawn from Model 1: in solid line under $P_0$, in dashed line under $P_1$.}
	\label{trajec}
\end{figure}

Table \ref{trajec} shows two trajectories drawn from Model 1, (in solid line, the trajectory drawn from $P_0$).

In both examples the training samples are of size 100 (50 trajectories drawn from each population). The outputs of the tables correspond to the average missclassification proportions (over 500 trajectories) of test samples of size 100 (50 generated from each population). 
The trajectories are discretized on a grid of 100 equispaced points. 

As for the choice of $k$, we have checked a reasonable range (according to the sample size) of values, in order to check ``robustness" with respect to $k$. We limit ourselves to odd values of $k$, from 3 to 9, just to avoid ties in the classifier output.

\begin{table}[ht]\label{tablesim1}
	\begin{center}
		
		\begin{tabular}{|c| c | c | c |}
			
			\hline 
			$k$    & $\mathbb{H}$ & $d_2$ & $d_\infty$ \\
			\hline 
			3  & .144  & .454   & .363 \\
			5  & .188  & .486  &  .450 \\
			7  & .222  & .496  &  .483 \\
			9  & .249  & .499  &  .494 \\
			\hline
			
		\end{tabular}
		\hspace{2cm}	\begin{tabular}{|c| c | c | c |}
			
			\hline 
			$k$    & $\mathbb{H}$ & $d_2$ & $d_\infty$ \\
			\hline 
			3  & .031  & .421  &  .275 \\
			5  & .046  & .473  &  .396 \\
			7  & .062  &  .491 &  .459 \\
			9  & .077  &  .497 &  .485\\
			\hline
			
		\end{tabular}
		\caption{Misclassification rates over 500 replications for different values of $k$ under Model 1 (left panel) and Model 2 (right).}
	\end{center}
\end{table}

The results are self-explanatory:  the classical distances have almost no discriminatory power in this example. The narrow noisy peaks are not suitable for them. This is in sharp contrast with the much better performance of the ${\mathbb H}$-distance.

\normalsize

\section{Real data examples}\label{real}

We will consider here two examples of binary classification based on functional data corresponding to mass spectra.  The ovarian cancer data is a bio-medical example. Hence the samples drawn from $P_0$ and $P_1$ correspond, respectively, to a control,``healthy''  group and to a ``patients group''; the aim is to assign a new coming individual with spectrum $x$ to one of these groups. The second example concerns food science: the goal is to investigate the capacity of mass spectra in order to discriminate between two varieties of coffee beans.

In both cases we have performed a similar experiment: the cross-validation (leave-one-out) proportions of correct classification have been computed for $k$-NN classifiers based on three different distances: the $L^2$-metric, $d_2$, the supremum metric $d_\infty$ and our Hausdorff-based distance ${\mathbb H}$.

The main goal of this study is just to check the possible usefulness of the ``visual'' distance ${\mathbb H}$ when compared with the classical choices $d_2$ and $d_\infty$. In principle, the idea was to handle the functional data themselves (or rather their available discretized versions) avoiding  the use of dimension reduction techniques via linear projections (principal components, partial least squares) or variable selection methods.In these examples, the available sample sizes are quite modest (especially in the second one). So the results must be interpreted with caution, just as useful hints of the performance of our proposal in real data problems. Obviously more research is needed.

\subsection{The ovarian cancer data}\label{ovarian}
These data correspond to mass spectra from blood samples of 216 women: 95 belong to the control group (CG) and the remaining 121 suffer from an ovarian cancer condition (OC). The use of mass spectra as a diagnostic tool in this situation is based on the fact that some proteins produced by cancer cells tend to be different (either in amount or in type) from that of the normal cells. These differences could be hopefully detected via mass spectrometry. We refer to  \cite{ban03} for a previous analysis of these data with a detailed discussion of their medical aspects. See also \cite{cues06} for further statistical analysis of these data. The raw data were defined on finite grids (of sizes varying between 320.000 and 360.000) on the interval $[700,12000]$. In order to facilitate the computational treatment we did some pre-processing: first,
we have restricted ourselves to the interval mass charge $[7000,9500]$, where  most peaks were concentrated. Second, we denoised the data by defining the spectra as $0$ at those points for which the value was smaller than 5 (this value was chosen after trying with several others). Third, in order to have all the spectra defined in a common equispaced grid (we took a grid of size 20001), we have smoothed them via a Nadaraya-Watson procedure. Finally, every function has been divided by its maximum, in order to have all the values scaled in the common interval $[0,1]$. This amounts to assume that the location of the peaks are important, but not the corresponding heights.

The results of our analysis are shown in Table \ref{tabesp} below. 

\begin{table}[h]
	\begin{center}
		\begin{tabular}{|c| c | c | c |}
			\hline 
			$k$    & $\mathbb{H}$ & $d_2$ & $d_\infty$ \\
			\hline 
			3  & .125  & .092   & .125  \\
			5  & .079  & .092   & .116  \\
			7  & .083  & .088   & .125  \\
			9  & .143  & .111   & .143  \\
			\hline
		\end{tabular}
		\caption{Classification error rates for the ovarian cancer data using $k$-NN classifiers based on three different distances.}
		\label{tabesp}
	\end{center}
\end{table}

We have kept the same values of $k$ considered in the simulation study. It can be seen that in all cases the "optimum" is either $k=5$ or $k=7$ and, for these choices, the Hausdorff based distance clearly outperforms the other two metrics $d_2$ and $d_\infty$.

\subsection{The coffee data}\label{coffee}

These data consist of 28 mass spectra (discretized in a grid of 286 values) corresponding to coffee beans of two varieties, \it Arabica\/ \rm and \it Robusta\rm. The respective sample sizes are 15 and 13. These data are available  from the web page
\url{http://www.cs.ucr.edu/~eamonn/time_series_data/}
of the University of California, Riverside. In this case the pre-processing consisted only on a rescaling of both axes to fit the data on $[0,1]^2$.

\begin{table}[h]
	\begin{center}
		\begin{tabular}{|c| c | c | c |}
			\hline 
			$k$    & $\mathbb{H}$ & $d_2$ & $d_\infty$ \\
			\hline 
			3  & .071  & .071   & .036  \\
			5  & .036  & .179    & 0  \\
			7  & .107  & .214   & .036  \\
			9  & .071 &  .25   &  .036 \\
			\hline
		\end{tabular}
		\caption{Classification error rates for the coffee data using $k$-NN classifiers based on three different distances.}
	\end{center}
\end{table}

In this case the ${\mathbb H}$-based classifiers are outperformed by those based on the supremum distance $d_\infty$ but are clearly better than those based on $d_2$. It is curious to note the unstable behavior of $d_2$, which is almost competitive in the cancer example but gets the worst performance both in the simulation study and in the coffee data  example. On the other hand, $d_\infty$ ranked clearly the last one in the cancer example. The ${\mathbb H}$-based methodology is never the worst one in the considered examples. Again, more detailed experiments are needed to confirm or refute this provisory findings.

\section{Concluding remarks}\label{conclu} The choice of a distance is particularly relevant when dealing with functional data. Not only some classifiers (as those of $k$-NN type or others based on depth measures) are defined in terms of distances, but also the theoretical properties (regarding consistency, convergence rates or asymptotic distributions) must be necessarily established in terms of a given distance in the sample space. Of course, in the finite-dimensional case, the use of different norms in the Euclidean sample space can lead to different results in a classification problem. However, the case for considering different norms in this finite-dimensional situation is not very strong, due to the well-known fact that all the norms are equivalent in finite-dimensional normed spaces.

Our generic suggestion here is to consider geometrically motivated distances in functional data. The specific proposal we make, ${\mathbb H}(f,g)$ is just one possibility; several other alternatives might be considered. The book by \cite{roc09} could suggest some ideas in this regard. While our theoretical and practical results with the distance ${\mathbb H}$ are encouraging, it is also quite clear that this metric suffers from some limitations: first, we are restricted to non-negative functions. Second, the extension of this idea to functions of several variables would probably involve considerable computational difficulties. Third, much more theoretical development is needed; in particular, the study of probability measures on the space $({\mathcal E},{\mathbb H})$ is essential if we want to use theoretical models combined with the distance ${\mathbb H}$. In fact, this need for a deeper mathematical development is a common feature for most chapters of the, still young, FDA theory.

\section*{Appendix}
\setcounter{section}{1}
\noindent \textbf{Proof of Proposition \ref{prophaus}}

To prove this Proposition we first must prove an auxiliary result:

\renewcommand{\thelema}{\textbf{A\thesection}}
\begin{lema} \label{lemhaus} If $f,g\in \mathcal{E}$, then there exist $u\in \partial H_f$ and $v\in \partial H_g$ such that
	\begin{equation} \label{hauslem1}
		\Haus(f,g)=d(u,H_g)=d(v,H_f)=\Vert u-v \Vert.
	\end{equation}
	\begin{proof}  We have, by definition of $\Haus$:
		$$\Haus(f,g)=d_H(H_f,H_g)=\max \bigg\{ \sup_{a\in H_f} d(a,H_g), \sup_{b\in H_g}d(b,H_f)\bigg\}.$$
		Assume ${\mathbb H}(f,g)>0$. Otherwise the result is trivial. Let us suppose by contradiction that there is no pair $(u,v)\in \partial
		H_f\times \partial H_g$ such that (\ref{hauslem1}) is fulfilled. In any case, the compactness of $H_f$ and $H_g$ guarantees the existence of $x=(x_1,x_2)$ and $y=(y_1,y_2)$  fulfilling (\ref{hauslem1}) but, according to our contradiction argument,
		either $x$ or $y$ must be
		an interior point. For example, if $x\in int(H_f)$, then $0<x_1<1$. We
		will see that $d\big((x_1,f(x_1)), H_g)\geq \Haus(f,g)$.  For every $t\in [0,1]$ such that  $|t-x_1|<\Haus(f,g)$
		let us denote $u^t=(t,u_2^t)$ and $v^t=(t,v^t_2)$ the
		intersection points  of $\partial B(x,\Haus(f,g))$ and the line
		$x_1=t$; with $u^t_2<v^t_2$. From the assumption on $x$,
		$d(x,H_g)={\mathbb H}(f,g)$.
		This entails that $\mathring{B}\big(x,\Haus(f,g)\big)\cap H_g=\emptyset$ and, since $H_g$ is a
		hypograph (which implies that if $(a,b)\in H_g$ then the segment
		joining $(a,b)$ and $(a,0)$ is included in $H_g$) it is clear
		that $g(t)\leq u^t_2$, for all $t\in [0,1]$ with $|t-x_1|<
		\Haus(f,g)$. Therefore, if we move upwards the point $x=(x_1,x_2)$ to $(x_1,f(x_1))$ (recall that from the USC assumption, $x_2\leq f(x_1)$), we have $\mathring{B}\big((x_1,f(x_1)),\Haus(f,g)\big)\cap
		H_g=\emptyset$ and then $d\big((x_1,f(x_1)), H_g\big)\geq
		\Haus(f,g)$.
		We cannot have
		$d\big((x_1,f(x_1)),
		H_g\big)>\Haus(f,g)$ since $(x_1,f(x_1))\in H_f$ and
		$\Haus(f,g)=d_H(H_f,H_g)$. So, we must have $d\big((x_1,f(x_1)), H_g\big)= \Haus(f,g)$ with $u:=(x_1,f(x_1))\in\partial H_f$. As a consequence, we must also have a point $v\in \partial H_g$ such that
		$\Vert u-v\Vert= \Haus(f,g)$. This contradicts the assumption we made about the non-existence of such a pair $(u,v)$.

	\end{proof}
\end{lema}

\noindent We now prove Proposition \ref{prophaus}.
\begin{proof}
	Let us denote
	\begin{footnotesize}
		$$d=\max \left\{\sup_{\left\{\substack{x=(x_1,x_2)\in \partial H_g \\ g(x_1)\geq f(x_1)}\right\}}d\big(x,\partial H_f\big),\ \sup_{\left\{\substack{y=(y_1,y_2)\in \partial H_f \\ f(y_1)> g(y_1)}\right\}} d\big(y,\partial H_g\big)\right\}.$$
	\end{footnotesize}
	The case $\Haus(f,g)=0$ is trivial, so let us assume
	$\Haus(f,g)>0$. We will first see that $\Haus(f,g)\leq d$. Since
	$H_f$ and $H_g$ are compact, there are two possibilities:
	\begin{itemize}
		\item[1)] there exists $x\in H_g$ such that $\Haus(f,g)=d(x,\partial H_f)$, or
		\item[2)] there exists $y\in H_f$ such that $\Haus(f,g)= d(y,\partial H_g)$.
	\end{itemize}
	
	Let us suppose that we are in the first case. By Lemma
	\ref{lemhaus} we can assume that $x=(x_1,x_2)\in \partial H_g$.
	Since $\Haus(f,g)>0$ and $H_f$ is a hypograph it must be $g(x_1)>
	f(x_1)$, then $x\in \big\{x=(x_1,x_2)\in \partial H_g:  g(x_1)\geq
	f(x_1)\big\}$ from where it follows that $\Haus(f,g)\leq d$. If we
	are in case 2, again by Lemma \ref{lemhaus}, we can assume
	$y=(y_1,y_2)\in
	\partial H_f$, as $\Haus(f,g)>0$ and $H_g$ is a hypograph it must
	be  $f(y_1)>g(y_1)$, then $y \in \big\{y=(y_1,y_2)\in \partial
	H_f: f(y_1)> g(y_1) \big\}$. from where it follows that
	$\Haus(f,g)\leq d$. The inequality ${\mathbb
		H}(f,g)\geq d$  follows directly from the definition of $\Haus$.
\end{proof}

We also use the following proposition in the algorithm to calculate ${\mathbb H}(f,g)$ for $f$ and $g$ continuous.

\renewcommand{\theprop}{\textbf{A\thesection}}
\begin{prop} \label{propdh} Let $f,g\in \mathcal{E}$ be continuous
	functions,
	let $u$ and $v$ be the points of Lemma \ref{lemhaus}. Then, there exist $t \in [0,1]$ and $s\in [0,1]$ such that $u=(t,f(t))$ and $v=(s,g(s))$.
	
	\begin{proof}
		Again, assume ${\mathbb H}(f,g)>0$. By Lemma
		\ref{lemhaus} $u\in \partial H_f$, $v\in \partial H_g$, and
		$\Haus(f,g)=\| u-v\|$. So it is enough to prove that
		$u=\big(t,f(t)\big)$ and $v=\big(s,g(s)\big)$. Since $f$ is
		continuous and $u\in \partial H_f$, there are four possibilities:
		(and the same holds for $v\in \partial H_g$) :
		\begin{itemize}
			\item[1.] $u$ is in the left border: $u=(0,u_2)$ with $u_2<f(0)$.
			\item[2.] $u$ is in the right border: $u=(1,u_2)$ with $u_2<f(1)$.
			\item[3.] $u$ is in the lower border: $u=(u_1,0)$ with $0\leq u_1\leq 1$.
			\item[4.] $u$ is in the upper border: $0\leq u_1\leq 1$ y $u_2=f(u_1)$.
		\end{itemize}
		We now prove that $u$ can only be in Case 4. It is clear that Case
		3 is not possible because both functions are non-negative. Cases 1
		and 2 are also excluded  following the ideas used in Lemma
		\ref{lemhaus}. For example, let us suppose that
		we are in Case 1 (see Figure \ref{hausdist3}).
		First observe
		that
		$\mathring{B}\big((0,f(0)),\Haus(f,g)\big)\cap
		H_g=\emptyset$; otherwise there would exist $(t_1,t_2)\in
		\mathring{B}\big((0,f(0)),\Haus(f,g)\big)\cap H_g$,
		then, the
		segment joining the points $(t_1,0)$ and $(t_1,t_2)$ (which is
		included in $H_g$) intersects
		$\mathring{B}\big(u,\Haus(f,g)\big)$. But this
		contradicts $\Haus(f,g)=d_H(H_f,H_g)$. So we
		conclude $d\big((0,f(0)),H_g\big)\geq \Haus(f,g)$.
		However $d\big((0,f(0)),H_g\big)> \Haus(f,g)$
		leads to a contradiction with the definition of $\Haus(f,g)$.
		Also, $d(u,H_g)= d\big((0,f(0)),H_g\big)$ leads to another
		contradiction. Indeed, if this were the case, we would have two
		points ($(0,u_2)$ and $(0,f(0))$) on the
		vertical 
		axis $x_1=0$ which are equidistant to the hypograph $H_g$. Then we
		have three possibilities:
		
		(a) $u_2<g(0)<f(0)$. This contradicts $d(u,H_g)=
		d\big((0,f(0)),H_g\big)$, since all the points $(0,u_3)$ with
		$u_2<u_3<g(0)$ belong to $H_g$.
		
		(b) $g(0)\leq u_2$: this contradicts the
		continuity of $g$ since $H_g$ must have a point in the boundary of
		$B\big((0,f(0)), d((0,f(0)),H_g)\big)$ and no point in the open
		ball $\mathring{B}\big((0,u_2), d\big((0,u_2),H_g\big)\big)$.
		
		(c) $g(0)\geq f(0)$: this is not compatible with
		$d(u,H_g)= d\big((0,f(0)),H_g\big)$.
		\begin{figure}[htbp]
			\centerline{\includegraphics[scale=.8]{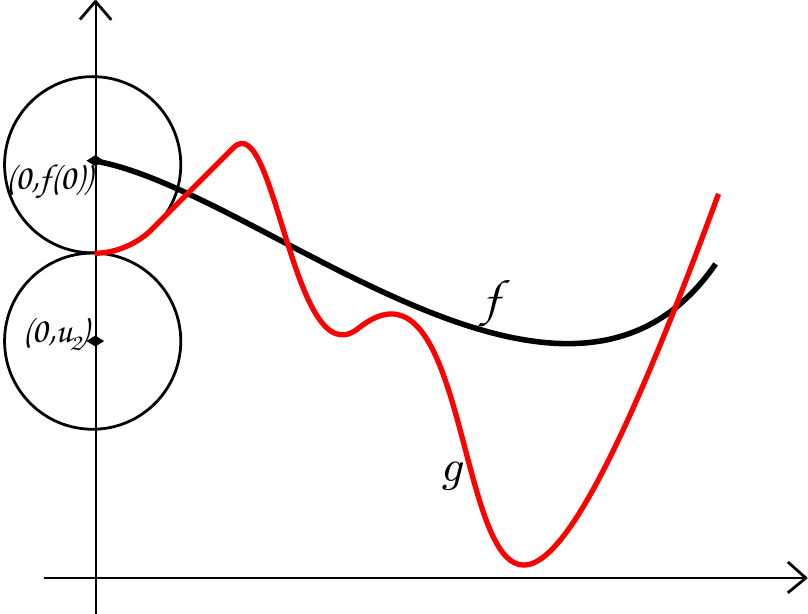}}
			\caption{We cannot have $d(u,H_g)= d\big((0,f(0)),H_g\big)$ with
				$u=(0,u_2)$ and $u_2<f(0)$.} \label{hausdist3}
		\end{figure}
	\end{proof}
\end{prop}

\noindent \bf Proof of Theorem \ref{topo}\rm.

(a) To state the local compactness we will in fact prove a slightly stronger property: we will show that any closed and bounded set in $({\mathcal E},{\mathbb H})$ is compact. Indeed, this would imply that the closed balls are compact. Since the family of balls with center at a given point is a local base, the local compactness will follow.

Since we are in a metric space compactness is equivalent to sequential compactness.  Let us take $\{f_n\}\subset \mathcal{E}$ a bounded sequence; we will prove that this sequence has necessarily a convergent subsequence. To see this, note that
the corresponding sequence of compact sets $H_{f_n}$ is  bounded. So it has convergent subsequence, which we may denote again by $H_{f_n}$, in the Hausdorff metric (since the closed and bounded sets are compact in the space of compact sets with the Hausdorff metric). Denote by $C$ the limit of that subsequence.
Therefore it is enough to prove that
\begin{eqnarray}
	&&\mbox{if }\big\{H_{f_n}\big\}_n \mbox{ is
		fulfils }
	H_{f_n}\rightarrow C \mbox{ for some compact set }C,\mbox{ then there }\nonumber\\
	&&\mbox{exists a USC function }
	f:[0,1]\rightarrow [0,\infty)\mbox{, such that }C=H_{f}.
	\label{HfC}
\end{eqnarray}

Let us take $(x,y)\in C$ and $(x_n,y_n)\in H_{f_n}$ converging to $(x,y)$; note that there exists at least one such sequence because $d_H(H_{f_n},C)\to 0$.  Now, since the $H_{f_n}$ are hypographs the vertical segment $[(x_n,0),(x_n,y_n)]$ joining the points $(x_n,0)$ and $(x_n,y_n)$ is included in $H_{f_n}$. So $H_{f_n}\rightarrow C$, implies
\begin{equation} \label{inc}
	\Big[\big(x,0\big),\big(x,\limsup f_n(x_n)\big)\Big]\subset C.
\end{equation}
Indeed, since $H_{f_n}$ is a hypograph,  $f_n(x_n)\geq y_n$. Then if we take $\limsup$ we obtain $y\leq \limsup f_n(x_n)$ and
$\big\{(x,z): 0\leq z \leq \limsup f(x_n)\big\}\subset C$.

Let us now define $f:[0,1]\rightarrow [0,\infty)$ by
$$
f(x)= \sup_{\{x_n\}:x_n\rightarrow x} \limsup f_n(x_n).
$$
Since $\{f_n\}$ is bounded, $f$ is well defined as a real-valued function. Let us prove that  $C=H_f$. Since $C$ is closed, we have, by (\ref{inc}),  $H_{f}\subset C$. Moreover, if $(x,y)\in C$, taking $(x_n,y_n)\rightarrow (x,y)$  with $(x_n,y_n)\in H_{f_n}$, $f(x)\geq \limsup f_n(x_n)\geq \limsup y_n=\, y$, we obtain $(x,y)\in H_f$.

It remains to prove that $f$ is USC. Suppose by contradiction that there exists $a \in [0,1]$ such that  $\limsup_{x\rightarrow a}f(x)>f(a)$. Then, we can  take a constant $\delta >0$ and a sequence $x_n\rightarrow a$, $x_n\neq a$ for all $n$  such that $f(x_n)>f(a)+\delta$ for $n$ large enough, say $n>n_0$. By the definition of $f$, for every $x_n$
we can take a sequence $z^k_n\rightarrow_k x_n$ (dependent on $n$), such that $f(x_n)=\lim_{z^k_n\rightarrow_k x_n}  \limsup f_k(z^k_n)$.

Given $\eps>0$, for every $n>n_0$ let us take take an increasing sequence $k(n)>n_0$ with
$$\big|z_n^{k(n)}-x_n\big|<\frac{1}{n}\ \mbox{and }\big|f_{k(n)}(z_n^{k(n)})-\limsup f_k(z^k_n)\big|<\eps,$$
that is,
$\big|f_{k(n)}(z_n^{k(n)})-f(x_n)\big|<\eps$.
But as $z_n^{k(n)}\rightarrow_n a$ this contradicts $f(x_n)>f(a)+\delta$ for $n>n_0$.

Completeness follows directly from the fact that the space of
compact sets endowed with the Hausdorff metric is complete,
together with (\ref{HfC}).

To prove separability,
let ${\mathcal P}_n$ be the set of all partitions of  $[0,1]$ defined by  $0=x_0<x_1<\dots<x_{n-1}<x_n=1$ where the $x_i$ are
rational numbers. Denote $\mathcal{P}=\cup_n {\mathcal P}_n$. Note that $\mathcal{P}$ in numerable.

Given a partition $\mathcal{P}\in
\mathcal{P}_n$ and a set $q_0,\dots,q_{n-1}$ of rational numbers, let us
define
\begin{equation} \label{dens}
	f_{\mathcal{P}}(x)=\begin{cases}
		q_0 & \text{ if } x \in [0,x_1)\\
		q_i &\text{ if } x\in (x_i,x_{i+1})\quad 1 \leq i \leq n-3\\
		q_{n-1} &\text{ if } x\in (x_{n-1},1]\\
		\max\{q_i,q_{i+1}\} &\text{ if } x=x_i \quad 1 \leq i \leq n-2
	\end{cases}
\end{equation}
It is immediately seen that this function  is USC and bounded. Let us see that the
(numerable)  set of all functions defined by \ref{dens}, for all possible partitions
$\mathcal{P}$ and rational values $q_i$ is dense in $\mathcal{E}$ with respect to
$\Haus$. Let $f$ be a  non-negative USC function  and take $\eps>0$. Consider
$\mathcal{P}\in \mathcal{P}_n$ a partition of the form $0=x_0<x_1<\dots<x_{n-1}<x_n=1$
where $x_i$ are rational numbers and such that
$\max_{i=0,\dots,n-1} |x_{i+1}-x_i|<\eps/2$. By Proposition
\ref{propusc} there exists $f_i=\max_{x\in [x_i,x_{i+1}]} f(x)$.
Let us take $q_0,\dots,q_{n-1}$ rational numbers  such that
$q_i>f_i$ and $q_i-f_i<\eps/2$ for all $i$. For this partition and
this set of rational numbers let us define $f_{\mathcal{P}}$ as in
(\ref{dens}). Now we claim that $\Haus(f_\mathcal{P},f)\leq\eps $. Indeed, it
is clear that $f_\mathcal{P}(x)>f(x)$ for all $x\in [0,1]$ so that
$H_f\subset H_{f_\mathcal{P}}$, and $\Haus(f_\mathcal{P},f)= \sup_{z\in H_{f_\mathcal{P}}}
d(z,H_f).$ Given $z=(z_1,z_2)\in H_{f_\mathcal{P}}$, there exists  $0\leq i_0
\leq n-1 $ such that $x_{i_0}\leq  z_1\leq x_{i_0+1}$. Now, choose
$t$ such that $x_{i_0}\leq  t\leq x_{i_0+1}$  and $f_{i_0}= f(t)$. We have
$z_2< f(t)+\eps/2$ and then $d(z,H_f)< \eps$.
Since $z$ was an arbitrary point in $H_{f_\mathcal{P}}$ we finally get
$\sup_{z\in H_{f_\mathcal{P}}} d(z,H_f)\leq\epsilon$.

\

(b) By Proposition \ref{propusc} (i) we know that there exists $z\in [0,1]$ such
that $f(z)=\max_{x\in [0,1]} f(x)$. As $\Haus(f_n,f)\rightarrow 0$ there
exist $x^n=(x_1^n,x_2^n)\in H_{f_n}$ such that $x_n\rightarrow (z,f(z))$.  Then,
$x^n_2\leq f_n(x_1^n)\leq \max_{x\in[0,1]} f_n(x)$ and, since $x_2^n\to f(z)$, we obtain
$$\max_{x\in [0,1]} f(x) = f(z)\leq \liminf_{n\rightarrow +\infty}  \max_{x\in[0,1]} f_n(x)
\leq \limsup_{n\rightarrow +\infty}  \max_{x\in [0,1]} f_n(x).$$
Finally, let us prove that  $\limsup_n \max_{x\in [0,1]} f_n(x) \leq \max_{x\in [0,1]} f(x)$.
Denote $z_0=\limsup _{n}  \max_{x\in [0,1]}
f_n(x)$. There exists $x_n\in[0,1]$ such that $f_n(x_n)\rightarrow
z_0$ with $f_n(x_n)= \max_{x\in [0,1]} f_n(x)$. Taking if necessary a
subsequence, we can assume that $x_n\rightarrow x_0\in[0,1]$. Since
$(x_0,z_0)\in H_f$ we have $f(x_0)\geq z_0$ then $\max_{x\in [0,1]}
f(x)\geq z_0$.

\

\noindent \bf Proof of Proposition \ref{ppp}\rm.

\begin{proof}
	This result is just a direct corollary from Th. 2 in \cite{cer06}
	(recall that the continuity of $\eta(x)$ is a sufficient condition
	for 
	(\ref{besicovitch2})), combined with the fact that the regression
	function $\eta_Q(x)={\mathbb P}(Y=1|X=x)$ (i.e., the regression
	function under $Q$) can be approximated by a continuous compact
	supported function; we use here the local compactness of
	${\mathcal E}$ (see \cite{fol99}, Proposition 7.9). Indeed, note
	that the joint distribution of $(X,Y)$ is completely determined by
	$\eta(x)={\mathbb P}(Y=1|X=x)$ and by the marginal distribution
	$\mu$ of $X$. Then, given $Q$, one can construct
	$P$ by just approximating $\eta_Q(x)={\mathbb P}(Y=1|X=x)$ by a continuous
	compact-supported function $\eta_P(x)$ which, without loss of
	generality, can be taken $0\leq \eta_P\leq 1$. Then, the
	distribution $P$ determined by $\eta_P$ and the
	marginal distribution $\mu$ of $X$
	is arbitrarily close to $Q$ (just
	taking $\eta_P$ close enough to $\eta$). Indeed,
	given any Borel set $C \subset \mathcal E \times \{0,1\}$,
	consider the sets $C_0 = \{x \in \mathcal E: (x,0) \in C\}$ and
	$C_1 = \{x \in \mathcal E: (x,1) \in C\}$. Then,
	$$
	Q(C) = \int _{C_0} (1-\eta_Q(x)) d\mu(x) + \int _{C_1}
	\eta_Q(x) d\mu(x),
	$$
	and
	$$
	P(C) = \int _{C_0} (1-\eta_P(x)) d\mu(x) + \int _{C_1}
	\eta_P(x) d\mu(x),
	$$
	which can be made arbitrarily close. 
\end{proof}

\end{document}